\documentstyle[11pt, pb-diagram]{article}
\input amssym.def
\newcommand\C{{\mbox{$\Bbb C$}}}
\newcommand\R{{\mbox{$\Bbb R$}}}

\newcommand\qed{\hfill $\Box$ }
\newtheorem{pr}{Proposition}

\newtheorem{th}{Theorem}
\newcommand{\be}{\begin{equation}}
\newcommand{\ee}{\end{equation}}
\newcommand{\bea}{\begin{eqnarray}}
\newcommand{\eea}{\end{eqnarray}}
\newcommand{\bean}{\begin{eqnarray*}}
\newcommand{\eean}{\end{eqnarray*}}
\newcommand\cD{{\cal D}}
\newcommand\cT{{\cal T}}
\newcommand\cE{{\cal E}}
\newcommand\cP{{\cal P}}

\newcommand\fa{{\frak a}}
\newcommand\fg{{\frak g}}
\newcommand\gt{\tau}
\newcommand\gk{\kappa}
\newcommand\var{\varepsilon}
\newcommand{\lu}{\underline{\lambda}}
\newcommand\mk{\medskip}
\title{Twistor spinors on Lorentzian symmetric spaces  }
\author{Helga Baum}
\date{\today}

\begin{document}

\maketitle
\begin{abstract}
We solve the twistor equation on all indecomposable Lorentzian symmetric spaces 
explicitly.
\end{abstract}

\section{Introduction}

Let $(M^n,g)$ be an oriented semi-Riemannian spin manifold with the
spinor bundle $S$. The {\it twistor operator} $\cD$ is defined as the
composition of the spinor derivative $\nabla^S$ with the projection
$p$ onto the kernel of the Clifford multiplication $\mu$

\[ \cD: \Gamma (S) \stackrel{\nabla^S}{\to} \Gamma (T^* M \otimes S )
\stackrel{g}{=} \Gamma (T M \otimes S) \stackrel{p}{\to} \Gamma
(\ker \mu)  . \]
The solutions of the conformally invariant equation $\,\cD \varphi = 0\,$ are 
called {\em twistor spinors}.
Twistor spinors were introduced by R.Penrose in General Relativity
(see \cite{Penrose:67}, \cite{Penrose/Rindler:86},
\cite{Nieuwenhuizen/Warner:84}). They are related to Killing vector
fields in semi-Riemannian supergeometry (see
\cite{Alekseevski/Cortes/ua:97}). In the last years essential results
concerning the geometry of Riemannian spin manifolds admitting twistor spinors 
were 
obtained by A.Lichnerowicz, T.Friedrich, K.Habermann, H-B.Rade\-macher, 
W.Kueh\-nel 
and other authors. For a survey on the literature cf. \cite{Friedrich:97}. In the 
Lorentzian setting there was established a relation between a special class of 
solutions of 
the twistor equation and the Fefferman spaces occuring in CR-geometry (cf. 
\cite{Lewandowski:91}, \cite{Baum:97}).\\
Let ${\cal T} (M^n,g)$ denote the space of all twistor spinors of
$(M,g)$. It is known that
\[ \dim {\cal T} (M^n,g) \le 2 \cdot 2^{\left[ \frac{n}{2} \right]} \]
(see \cite{Baum/Friedrich/ua:91}). If $(M^n,g)$ is conformally flat
and simply connected, then one has  $\,\dim {\cal T} (M^n,g)= 2
\cdot 2^{\left[ \frac{n}{2} \right]}\,$. In the present paper we
determine the twistor spinors on all indecomposable Lorentzian symmetric 
spaces expilicity. In particular, we prove:
\begin{enumerate}
\item
If $(M^n,g)$ is an indecomposable non-conformally flat
Lorentzian symmetric spin manifold of dimension $n \ge 3$, then each
twistor spinor is parallel and
$\,\dim {\cal T} (M^n,g) = q \cdot 2^{\left[ \frac{n}{2} \right]}\,$, where $q= 
\frac{1}{2}, \frac{1}{4}$ or 0,
depending on the
fundamental group $\pi_1 (M)$ and on the spin structure of
$(M^n,g)$.
\item
If $(M^n,g)$ is an indecomposable conformally flat Lorentzian
symmetric spin manifold of dimension $n \ge 3$ and non-constant
sectional curvature, then
$\,\dim {\cal T} (M^n,g)= q \cdot 2^{\left[
\frac{n}{2} \right]}\,$, where $q = 2, \frac{3}{2}, 1, \frac{3}{4}$ or
$0$ , depending on $\pi_1 (M)$ and on the spin structure.
\item
If $(M^n,g)$ is a Lorentzian symmetric spin manifold of
dimension $n \ge 3$ and constant sectional curvature, then $\,\dim
{\cal T} (M^n,g) =q \cdot 2^{\left[ \frac{n}{2} \right]}\,$, where
$q=2,1,$ or 0, depending on $\pi_1 (M)$ and on the spin structure.
\end{enumerate}
The calculations are based on the fact, that the indecomposable
Lorentzian symmetric spaces are completely classified (see
\cite{Cahen/Wallach:70}). \\

\section{Lorentzian symmetric spaces}

Let us first recall the description of Lorentzian symmetric spaces.
A connected semi-Riemannian manifold $(M,g)$ is called {\em
indecomposable}, if there is no proper, nondegenerate
subspace of $T_x M$ invariant under the action of the holonomy group
$\mbox{Hol}_x (g)\,$. Each simply connected semi-Riemannian symmetric space is 
isometric to a product $\,M_0 \times M_1 \times \dots \times M_r\,$,
where $M_i, \,i=1, \dots
, r\,$, are indecomposable simply connected semi-Riemannian symmetric
spaces of
dimension $\geq 2$ and $M_0$ is semi-Euclidean.\\
Let $(M^n,g)$ be a Lorentzian symmetric space. By $G(M)$ we denote
the group of transvections of $(M^n,g)$ and by $\fg$ its
Lie algebra. One has the following structure result:

\begin{th} {\em (\cite{Cahen/Wallach:70})}\\
Let  $(M^n,g)$ be an indecomposable Lorentzian symmetric space of
dimension $n \ge 2$. Then the Lie algebra $\fg$ of the
transvection group of $(M^n,g)$ is eigther semi-simple or
solvable.
\end{th}

Let $\lu =( \lambda_1, \ldots, \lambda_{n-2})$ be
an $(n-2)$-tupel of real numbers $\lambda_j \in {\Bbb R} \backslash
\{0\}$ and let us denote by $M^n_{\lu}$ the
Lorentzian space $M^n_{\lu}:= ( {\Bbb R}^n,
g_{\lu} )$, where
\[
 (g_{\lu})_{(s,t,x)} := 2ds\,dt +
\sum\limits^{n-2}_{j=1} \lambda_j x^2_j \,ds^2 +
\sum\limits^{n-2}_{j=1} dx_j^2 .
 \]
If $\lu_{\pi} =( \lambda_{\pi(1)}, \ldots ,
\lambda_{\pi (n-2)} )$ is a permutation of $\lu$ and
$c >0$, then $M^n_{\lu}$ is isometric to $M^n_{c\lu_{\pi}}$.

\begin{th} {\em (\cite{Cahen/Wallach:70}, \cite{Cahen/Parker:80})}\\
Let $(M^n,g)$ be an indecomposable solvable Lorentzian symmetric
space of dimension $n \ge 3$. Then $(M^n,g)$ is isometric to
$\,M^n_{\lu}/A\,$, where $\,\lu \in ({\Bbb R}
\backslash \{ 0 \})^{n-2}\,$ and $A$ is a discrete subgroup of the
centralizer $Z_{I(M_{\lu})}
(G(M_{\lu}))$ of the transvection group
$G(M_{\lu})$ in the isometry group
$I(M_{\lu})$ of $M^n_{\lu}$.
\end{th}

For the centralizer
$\,Z_{\lu}:=Z_{I(M_{\lu})}
(G(M_{\lu}))\,$ it is known:

\begin{th} {\em  (\cite{Cahen/Kerbrat:78})}\\
Let $\lu =( \lambda_1, \ldots, \lambda_{n-2})$ be a tupel of
non-zero real numbers.
\begin{enumerate}
\item
If there is a positive $\lambda_i$ or if there are two numbers $\lambda_i, 
\lambda_j$ such that $\frac{\lambda_i}{\lambda_j} \not\in {\Bbb
Q}^2$, then $\,Z_{\lu} \simeq {\Bbb R}\,$ and $\varphi
\in Z_{\lu}\,$ if and only if $\varphi (s,t,x)=(s,t+ \alpha,
x)$, $\alpha \in {\Bbb R}$.
\item
Let $\lambda_i =- k^2_i <0$ and $\,\frac{k_i}{k_j} \in {\Bbb
Q}\,$ for all $i,j \in \{ 1, \ldots , n-2\}$. Then $\,\varphi \in
Z_{\lu}\,$ if and only if
\[
\varphi (s,t,x) =(s+ \beta, t+ \alpha, (-1)^{m_1} x_1, \ldots ,
(-1)^{m_{n-2}} x_{n-2}) ,
 \]
where $\alpha \in {\Bbb R},\, m_1, \ldots, m_{n-2} \in {\Bbb Z}\,$ and
$\,\beta = \frac{m_i \cdot \pi}{k_i}\,$ for all $i=1, \ldots, n-2 .$
\end{enumerate}
\end{th}
\vspace{0.5cm}
Let us denote by $S^n_1 (r)$ the pseudo-sphere
\[
S^n_1 (r) := \left\{ x \in {\Bbb R}^{n+1,1} \mid \langle
x,x \rangle_{n+1,1} = -x^2_1 + x^2_2 + \ldots + x^2_{n+1} =r^2 \right\} 
\,\subset {\Bbb R}^{n+1,1}
\]
and by $H^n_1 (r)$ the pseudo-hyperbolic space
\[
 H^n_1 (r):= \left\{ x \in {\Bbb R}^{n+1,2} \mid \langle x,x
\rangle_{n+1,2} = -x^2_1 - x^2_2 + x_3^2 + \ldots + x^2_{n+1} = -r^2
\right\} \subset {\Bbb R}^{n+1,2}
 \]
with the Lorentzian metrics induced by $\langle \cdot , \cdot \rangle_{n+1,1}$ 
and $\langle  \cdot , \cdot \rangle_{n+1,2}$, respectively.

\begin{th} {\em  (\cite{Cahen/Leroy/ua:90}, \cite{Wolf:84})}\\
Let $(M^n,g)$ be an indecomposable semi-simple Lorentzian symmetric
space of dimension $n \ge 3$. Then $(M^n,g)$ has constant sectional
curvature $k \not= 0$. Therefore, it is isometric to $S^n_1 (r)
/_{\{ \pm I \}}$ or $S^n_1 (r)\,,\, (k = \frac{1}{r^2} >0)$, or to
a Lorentzian covering of $H^n_1 (r)/_{\{\pm I \}} \,,\, (k= -
\frac{1}{r^2} <0)$.
\end{th}

\section{Spinor representation}

For  concrete calculations we will use the following
realization of the spinor representation.
Let $\,\mbox{Cliff}_{n,k}\,$ be the Clifford algebra of
$\,(\R^n,-\langle \cdot, \cdot \rangle_k)\,$, where $\,\langle \cdot, \cdot 
\rangle_k\,$ is the scalar product
$\,\, \langle
x,y\rangle_k:=-x_1y_1-\ldots-x_ky_k+x_{k+1}y_{k+1}+\ldots+x_ny_n\,\,$.
For the canonical basis
$(e_1,\ldots,e_n)$ of $\R^n$ one has the following relations in
$\,\mbox{Cliff}_{n,k}\,:\,\,
 e_i\cdot e_j+e_j\cdot
e_i=-2\var_j\delta_{ij},\,\,$  where $ \var_j=\left\{\begin{array}{rl}
-1 & j\le k\\ 1 & j> k \,\end{array}\right. $.
Denote $\,\,\tau_j = \left\{\begin{array}{ll} i & j\le k\\ 1 &
j>k\end{array}\right.\,\,$ and
\begin{displaymath}
U = \left(\begin{array}{cc} i & 0\\ 0 & -i\end{array}\right),\quad
V=\left(\begin{array}{cc} 0 & i\\ i & 0\end{array}\right),\quad
E=\left(\begin{array}{cc} 1 & 0\\ 0 & 1\end{array}\right),\quad
T=\left(\begin{array}{cc} 0 & -i\\ i &
0\end{array}\right).
\end{displaymath}
If $n=2m$ is even, we have an isomorphism
\[
 \phi_{2m,k}\,:\,\mbox{Cliff}^{\Bbb C}_{2m,k}
\stackrel{\sim}{\longrightarrow} M(2^m;\C)
 \]
given by the Kronecker product
\begin{eqnarray}\label{phi1}
\begin{array}{llll} \phi_{2m,k}(e_{2j-1}) & = &
\tau_{2j-1} &E\otimes\ldots\otimes E \otimes U \otimes
T\otimes\ldots\otimes T\\ \phi_{2m,k}(e_{2j}) & =
& \tau_{2j} &E\otimes\ldots\otimes E\otimes V \otimes
\underbrace{T\otimes\ldots\otimes T}_{j-1} \end{array}.
\end{eqnarray}
If $n=2m+1$ is odd and $k<n$, we have the isomorphism
\[
\phi_{2m+1,k}\,:\, \mbox{Cliff}^{{\Bbb C}}_{2m+1,k}
\stackrel{\sim}{\longrightarrow} \,M(2^m;\C)
\oplus M(2^m;\C)
\]
given by
\begin{eqnarray}\label{phi2}
\begin{array}{llll}
\Phi_{2m+1,k}(e_j) & = & (\Phi_{2m,k}(e_j)\,,\, \Phi_{2m,k}(e_j)), \;\;\; 
j=1,\ldots, 2m\\[2ex]
\Phi_{2m+1,k}(e_n) & = & \gt_n (i\,\, T\otimes\cdots \otimes T\,,
\,-i\,\,T\otimes\cdots \otimes T).
\end{array}
\end{eqnarray}
Let
$\mbox{Spin}(n,k)\subset \mbox{Cliff}_{n,k}$ be the spin group. The spinor 
representation is given
by \[ \gk_{n,k}=\hat{\phi}_{n,k}|_{\mbox{Spin}(n,k)}:\mbox{Spin}(n,k)
\longrightarrow \mbox{GL} (\C^{2^{\left[\frac n2\right]}}), \]
where $\,\hat{\phi}_{2m,k} = \phi_{2m,k}\,$ and $\,\hat{\phi}_{2m+1,k} = pr_1 
\circ \phi_{2m+1,k}\,$. 
We denote this representation by
$\,\Delta_{n,k}\,$. If $n=2m$, $\,\Delta_{2m,k}\,$ splits into the sum
$\,\Delta_{2m,k}=\Delta^+_{2m,k}\oplus \Delta^-_{2m,k}\,$, where
$\,\Delta^\pm_{2m,k}\,$ are the eigenspaces of the endomorphism
$\,\phi_{2m,k}(e_1\cdot\ldots\cdot e_{2m})\,$ to the eigenvalue $\pm
i^{m+k}$.
Let us denote by $u(\delta)\in\C^2$ the vector $\,
u(\delta)=\frac{1}{\sqrt{2}}{1\choose -\delta i},\,\,\delta=\pm 1
\,\,$, and let
 \be\label{u}
u(\delta_1,\ldots,\delta_m)=u(\delta_1)\otimes\ldots\otimes
u(\delta_m)\qquad \delta_j=\pm 1. \ee Then
$\;(u(\delta_1,\ldots,\delta_m)\,|\,\prod\limits^m_{j=1}\delta_j=\pm
1)\;$ is an orthonormal basis of $\,\Delta^\pm_{2m,k}\,$ with respect
to the standard scalar product of $\C^{2^m}$.\\

\section{General properties of twistor spinors}

In this section we recall some properties of twistor spinors which
we will need in the following calculations. For proofs see
\cite{Baum/Friedrich/ua:91}. Let $(M^n,g)$ be an oriented
semi-Riemannian spin manifold of dimension $n \ge 3$. We denote by
$S$ the spinor bundle of $(M^n,g)$, by $\nabla^S : \Gamma (S) \to
\Gamma (T^* M \otimes S)$ the spinor derivative given by the
Levi-Civita connection of  $(M^n,g)$ and by $D: \Gamma (S) \to
\Gamma (S)$ the Dirac operator on $S$. Let $p:TM \otimes S \to TM
\otimes S$ be the projection onto the kernel of the Clifford
multiplication $\mu$. The {\em twistor operator}
\[
\cD: \Gamma (S) \stackrel{\nabla^S}{\to} \Gamma (T^*M \otimes S)
\stackrel{g}{=} \Gamma  (TM \otimes S) \stackrel{p}{\to} \Gamma
(\ker \mu)
 \]
is locally given by
\[
\cD \varphi = \sum\limits^n_{j=1} \var_j s_j \otimes (\nabla^S_{s_j}\varphi 
+ \frac{1}{n} s_j \cdot D \varphi) ,
\]
where $(s_1, \ldots, s_n)$ is a local orthonormal basis and $\var_j
=g(s_j, s_j)= \pm 1$. A spinor field $\varphi \in \Gamma (S)$ is
called {\em twistor spinor} if $\,\cD \varphi =0$.
\mk
\begin{pr}
For a spinor field  $\varphi\in
\Gamma(S)$ the following conditions are equivalent:\\[0.1cm]
\begin{tabular}{cl}
1. & $\varphi$ is a twistor spinor.\\[0.1cm]
2. & $\varphi$ satisfies
the so-called twistor equation
\end{tabular}
\vspace{-0.05cm}\\
\begin{equation} \label{4.1}
\nabla^S_X\varphi+\frac{1}{n}X\cdot D\varphi=0
\end{equation}
\vspace{-0.4cm}\\
\begin{tabular}{cl}
   & for all vector fields $X$.\\[0.1cm]
3. & There exists a spinor field $\psi\in\Gamma(S)$ such that
\end{tabular}
\vspace{-0.05cm}\\
\be\label{4.2}
\psi=g(X,X)X\cdot\nabla^S_X\varphi \ee
\vspace{-0.4cm}\\
\begin{tabular}{cl}
   & $\;\;\;$for all vector fields $X$ with $\,|g(X,X)|=1$.
\end{tabular}
\end{pr}

The dimension of the space ${\cal T}(M^n,g)$ of all twistor spinors
is conformally invariant and bounded by
\[
\dim {\cal T}(M^n,g) \le  2 \cdot 2^{\left[ \frac{n}{2} \right] } .
\]
If $(M^n,g)$ is simply connected and conformally flat, $\,\dim {\cal
T}(M^n,g)= 2\cdot 2^{\left[ \frac{n}{2} \right] }\,$. \\
In particular, the
twistor spinors on the semi-Euclidean space $\, \R^{n,k}:= ({\Bbb R}^n, 
\langle  \cdot, \cdot \rangle_{n,k})\,$ are the functions
\[
\cT(\R^{n,k}) = \left\{ \varphi \in C^{\infty} ({\Bbb
R}^{n,k} , \Delta_{n,k}) \,\mid \,\varphi (x) = u+x \cdot v\,; \,\,u,v \in 
\Delta_{n,k} \right\}.
\]
Let $R$ be the scalar curvature and Ric the Ricci curvature of
$(M^n,g)$. \linebreak $K:TM \to TM$ denotes the $(1,1)$-Schouten tensor of 
$(M^n,g)$
\[
K(X) = \frac{1}{n-2} \left( \frac{R}{2(n-1)} X - \mbox{Ric} (X)
\right) .
 \]
Furthermore, let $W$ be the $(4,0)$-Weyl tensor of $(M,g)$ and let us
denote by the same symbol the corresponding $(2,2)$-tensor field $W
: \Lambda^2 M \to \Lambda^2 M$. Then we have
\mk
\begin{pr}
Let $\varphi \in \Gamma (S)$ be a twistor spinor. Then
\begin{eqnarray}
\label{4.4} \nabla^S_X D
\varphi & = & \frac{n}{2} K(X) \cdot \varphi\;,\\
\label{4.5} W (\eta )
\cdot \varphi &=& 0 \;,
\eea
for all vector fields $X$ and 2-forms $\eta$.
\end{pr}
\mk

Finally, we recall two possibilities to obtain new manifolds with
twistor spinors from a given one.\\
Let $(\tilde{M}^n, \tilde{g})$ be a simply connected
parallelizable semi-Riemannian manifold and let $A \subset I(\tilde{M}, 
\tilde{g})$ denote a discrete subgroup of orientation preserving isometries 
of $(\tilde{M}, \tilde{g})$. We trivialize the spin structure of $(\tilde{M}, 
\tilde{g})$ with respect to a fixed global orthonormal basis field $\fa=(\fa_1, 
\dots , 
\fa_n)$.
For $\gamma \in A$
we denote by $\Gamma (x) \in SO(n,k)$ the matrix of $d \gamma_x$
with respect to $\fa (x)$ and $\fa (\gamma (x))$. Then there are
two lifts $\tilde{\Gamma}^{\pm}$ of $\Gamma$ into $Spin(n,k)$
\[
\begin{diagram}
\node{} \node{Spin (n,k)} \arrow{s,b}{\lambda} \\
\node{\tilde{M}} \arrow{ne,t}{\tilde{\Gamma}^{\pm}}
\arrow{e,b}{\Gamma} \node{SO(n,k)}
\end{diagram}
\]
Let $\cE (A)$ be the set of all left actions of $A$ on $\tilde{M}
\times Spin (n,k)$ such that
\[
\varepsilon (\gamma)(x, a)=( \gamma (x), \varepsilon
(\gamma, x) \cdot a)  \qquad \mbox{and} \qquad
 \varepsilon (\gamma , x)= \tilde{\Gamma} (x)^{\pm} .
 \]
This set of left actions $\cE (A)$ corresponds to the set of spin
structures of the oriented semi-Riemannian manifold $M=
\tilde{M}/A$. The spinor bundle on $M$ corresponding to the spin structure 
$\,\varepsilon \in \cE (A)\,$ is given by
\[
 S_{\varepsilon} = \left( \tilde{M} \times
\Delta_{n,k} \right) \Big/_{\textstyle{\varepsilon}}
 \]
where $\,\varepsilon (\gamma) (x,v) =(\gamma (x),
\varepsilon (\gamma, x) \cdot v)\,$ for all $\gamma \in
A , \, (x,v) \in \tilde{M} \times
\Delta_{n,k}\, $. Hence, the spinor fields on $M$ corresponding to $\varepsilon 
\in \cE (A)$ are given by the $\varepsilon$-invariant functions
\begin{eqnarray*}
\Gamma (S_{\varepsilon} ) &=& C^{\infty} (\tilde{M}, \Delta_{n,k}
)^{\textstyle{\varepsilon}} :=
\left\{ \varphi \in C^{\infty} (\tilde{M}, \Delta_{n,k})\,\mid\,
\varphi (\gamma (x))= \varepsilon (\gamma,x) \cdot
\varphi (x)\right\},
\end{eqnarray*}
and for the twistor spinors on $M$ with the spin structure
$\varepsilon$ we have
\mk
\begin{pr}
The twistor spinors on $\,M = \tilde{M}/A\,$ with respect to the
spin structure $\varepsilon \in \cE (A)$ are given by
\[
{\cal T}((M,g) , \varepsilon) = \{ \varphi \in {\cal T}(\tilde{M},
\tilde{g}) \,\mid \,\, \varphi \, \, \mbox{ist
$\varepsilon$-invariant} \,\,\}.
 \]
\end{pr}

\mk
Let $(M^{n+1}, g)$ be a semi-Riemannian spin manifold with spinor bundle $S_M$ and 
let $F^n \subset M^{n+1}$ be a
non-degenerate oriented hypersurface in $M^{n+1}$. We denote by
$\eta :F \to TM$ the Gauss map of $F, \; \,\kappa (\eta) := g( \eta,
\eta)= \pm 1$. It is well known, that the spinor bundle $S_F$ of
$(F, g|_{F})$ with respect to the spin structure on $F$ induced by
the embedding, is  isomorphic to $S_M |_F$ in case of even
dimension $n$ and to $S_M^{\pm} |_F$ in case of odd dimension $n$.
Using this identification the Clifford multiplication and the spinor
derivative are expressed by
\begin{eqnarray*}
X \cdot (\varphi|_F) &=&(X \cdot \varphi)|_F\\
\nabla_X^{S_F} (\varphi|_F)&=& \left( \nabla_X^{S^{(\pm)}_M} \varphi
+ \frac{1}{2} \kappa (\eta)\, \nabla^M_X \eta \cdot \eta \cdot
\varphi \right) \Big|_F,
\end{eqnarray*}
where $\,\varphi \in \Gamma (S^{(\pm)}_M) , X \in TF$ and $\varphi|_F$ always 
means the spinor field in $\Gamma(S_F)$ corresponding to $\varphi$ with respect 
to the 
above mentioned isomorphism.
\mk

\begin{pr}
If $F^n \subset M^{n+1}$ is an umbilic hypersurface and $\varphi
\in \Gamma (S^{(\pm)}_M)$ is a twistor spinor on $M$, then
$\varphi|_F \in \Gamma (S_F)$ is a twistor on $F$.
\end{pr}

\underline{Proof:} Let $\lambda \in C^{\infty} (F)$ be the function
satisfying $\,\nabla^M_X \eta = \lambda \, X\,$. Then
\[
 \nabla_X^{S_F} (\varphi |_F)= \left( \nabla_X^{S^{(\pm)}_M}
\varphi + \frac{1}{2} \kappa (\eta) \lambda X \cdot \eta \cdot
\varphi \right) \Big|_F .
\]
If $\varphi \in \Gamma (S^{(\pm)}_M)$ is a twistor spinor, from
Proposition 1, (\ref{4.1}), follows
\[
\kappa (X) \, X \cdot \nabla_X^{S^{(\pm)}_M} \varphi =
\frac{1}{n+1} D_M^{(\pm)}\varphi
\]
for each $X \in TF$ with $\kappa (X)= g(X,X)= \pm 1$. Hence,
\[
\kappa (X) \, X \cdot \nabla_X^{S_F} (\varphi |_F )= \left(
\frac{1}{n+1} D^{(\pm)}_M \varphi - \frac{1}{2} \kappa (\eta)
\lambda \,\eta \cdot \varphi \right) \Big|_F .
\]
The right hand side is independend of $X \in TF$. Therefore,
according to Proposition 1, (\ref{4.2}), $\varphi |_F$ is a twistor spinor on 
$F$. \qed\\

\section{Twistor spinors on indecomposable, non-conformally flat
Lorentzian symmetric spaces}

Let us first consider the simply connected solvable Lorentzian
symmetric space $M^n_{\lu}=({\Bbb R}^n, g_{\lu})$, where
\[
 (g_{\lu})_{(s,t,x)} = 2\,ds\,dt + \sum\limits^{n-2}_{j=1} \lambda_j
x_j^2\, ds^2 + \sum\limits^{n-2}_{j=1} dx_j^2 , \]

and $\lu = (\lambda_1, \ldots , \lambda_{n-2}), \, \, \lambda_i \in
{\Bbb R} \backslash \{ 0 \}, \, \, n \ge 3$. Let $\Lambda_0 := -
\sum\limits^{n-2}_{j=1} \lambda_j$. We fix the following global
orthonormal basis on $M^n_{\lu}$:
\begin{eqnarray*}
\fa_{\overline{0}} (y)&:=& \textstyle{\frac{\partial}{\partial s} (y) - 
\frac{1}{2} \left( \sum\limits^{n-2}_{j=1} \lambda_j x_j^2 +1
\right) \frac{\partial}{\partial t} (y)} \\
\fa_0 (y) &:=& \textstyle{\frac{\partial}{\partial s} (y) - \frac{1}{2} 
\left( \sum\limits^{n-2}_{j=1} \lambda_j x_j^2 - 1 \right)
\frac{\partial}{\partial t} (y)}\\
\fa_j (y) &:=& \textstyle{\frac{\partial}{\partial x_j} (y) \hspace{3cm} 
j=1, \ldots, n-2} ,
\end{eqnarray*}
where $y =(s,t,x_1, \ldots , x_{n-2} ) \in M^n_{\lu}$. The vector
field $\,V(y):= \frac{\partial}{\partial t} (y)\, $ is isotropic and
parallel. The Ricci tensor of $M^n_{\lu}$ is given by
\[
 \mbox{Ric} \, (X)= \Lambda_0 \cdot g_{\lu}(X,V)V ,
\]
the scalar curvature $R$ vanishes. Therefore, the Schouten tensor
satisfies
\bea \label{5.1}
K(X)= - \textstyle{\frac{1}{n-2}}\, \Lambda_0 \cdot g_{\lu} (X,V)\,V .
\eea
For the Weyl tensor $\,W: \Lambda^2 M_{\lu} \to \Lambda^2 M_{\lu}\,$ one has 
\bea
W( \fa_{\overline{0}} \wedge \fa_j) & = & W( \fa_0
\wedge \fa_j) \,= \, (\lambda_j + \textstyle{\frac{1}{n-2}} \Lambda_0)\, V 
\wedge
\fa_j , \qquad \, j= 1, \ldots, n-2 \quad \label{5.2}\\
W( \fa_{\alpha} \wedge \fa_{\beta}) & = & 0 \hspace{5.5cm}\mbox{for all 
other indices}\, \,\alpha, \beta, \nonumber
\eea
where $TM_{\lu}$ is identified with $T^*M_{\lu}$ using the metric
$g_{\lu}$. In particular, $M_{\lu}$ is conformally flat if and only
if $\lu =(\lambda, \ldots , \lambda)$, $\lambda \in {\Bbb R}
\backslash \{ 0 \}$.\\
Since $M_{\lu}$ is simply connected, it has an uniquely determined
spin structure. We trivialize this spin structure using the
global orthonormal basis $(\fa_{\overline{0}}, \fa_0, \fa_1,
\ldots , \fa_{n-2})$ and identify the spinor fields with the
smooth functions $C^{\infty} (M_{\lu}, \Delta_{n,1})$. The spinor
derivative is defined by

\[ \nabla^S_X \varphi =X(\varphi) + \frac{1}{2} \sum\limits_{1 \le
k < l \le n} \var_k \var_l \,g (\nabla^{LC}_X s_k,s_l)\, s_k \cdot
s_l \cdot \varphi , \]

where $(s_1, \ldots , s_n)$ is a local orthonormal basis and
$\var_j =g (s_j, s_j) = \pm 1$. This gives for the spinor derivative on $M_{\lu}$ 
\bea
\nabla^S_{\frac{\partial}{\partial t}} \varphi &  =  &
\textstyle{\frac{\partial}{\partial t} \varphi} \label{5.3}\\
\nabla^S_{\frac{\partial}{\partial s}} \varphi & = &
\textstyle{\frac{\partial}{\partial s} \varphi + \frac{1}{2}
\sum\limits^{n-2}_{j=1} \lambda_j x_j \fa_j \cdot V \cdot \varphi}
\label{5.4}\\
\nabla^S_{\frac{\partial}{\partial x_j}} \varphi & = &
\textstyle{\frac{\partial}{\partial x_j} \varphi}  \label{5.5}\\
\nabla^S_{\fa_{\overline{0}}} \varphi & = & \textstyle{\fa_{\overline{0}} 
(\varphi) + \frac{1}{2} \sum\limits^{n-2}_{j=1} \lambda_j x_j
\fa_j \cdot V \cdot \varphi}  \label{5.6}\\
\nabla^S_{\fa_0} \varphi & = & \textstyle{\fa_0 (\varphi) + \frac{1}{2} 
\sum\limits^{n-2}_{j=1} \lambda_j x_j \fa_j \cdot V \cdot \varphi}.
\label{5.7}
\eea
The vector space $\Delta_{n,1}$ is isomorphic to $\,\Delta_{n-2,0}
\otimes {\Bbb C}^2\,$. Let us denote by $\Delta_V$ the subspace
\[
 \Delta_V := \Delta_{n-2,0} \otimes {\Bbb C} u(-1) \quad \subset
\Delta_{n,1}  .
\]
Using the formulas (\ref{phi1}), (\ref{phi2}) and (\ref{u}) one obtains that a 
spinor field $\varphi
\in C^{\infty} (M_{\lu}, \Delta_{n,1})$ satisfies $\,V \cdot \varphi
=0\,$ if and only if the image of $\varphi$ lies in $\Delta_V$.
\mk
\begin{pr}
The space $\cP(M_{\lu})$ of parallel spinors of $M_{\lu}$ is
\begin{eqnarray*}
\cP (M_{\lu})= \left\{ \varphi \in C^{\infty} (M_{\lu},
\Delta_{n,1} ) \mid \,
 \varphi = \mbox{constant} \, \in \Delta_V \right\}.
\end{eqnarray*}
In particular, $\,\dim \cP(M_{\lu})= \textstyle{\frac{1}{2}} \cdot
2^{\left[\frac{n}{2} \right]
}$.
\end{pr}
\underline{Proof:} From (\ref{5.3}) - (\ref{5.5}) it follows that $\varphi \in 
C^{\infty} (M_{\lu}, \Delta_{n,1})$ is parallel if and only if
$\varphi$ depends only on $s$ and
\bea
\textstyle{\frac{\partial \varphi}{\partial s} = -
\frac{1}{2} \sum\limits^{n-2}_{j=1} \lambda_j x_j\, \fa_j \cdot V
\cdot \varphi}.  \label{5.8}
\eea
Therefore,
\[ \textstyle{0= \frac{\partial^2 \varphi}{\partial x_k \partial s} = - 
\frac{1}{2} \lambda_k \fa_k \cdot V \cdot \varphi} . \]
Since $\lambda_k \not= 0$ and $\fa_k$ is space-like, this yields
$V \cdot \varphi =0$. Hence, because of (\ref{5.8}), $\varphi$ has to be 
constant.
\qed
\mk
\begin{pr}
Let $M_{\lu}$ be non-conformally flat. Then each twistor spinor on
$M_{\lu}$ is parallel. In particular,
\[
\dim {\cal T} (M_{\lu})= \textstyle{\frac12}\cdot 2^{\left[ \frac{n}{2} 
\right]} . \]
\end{pr}
\underline{Proof:} Let $\varphi \in C^{\infty} (M_{\lu},
\Delta_{n,1})$ be a twistor spinor. Then according to (\ref{4.5}) of
Proposition 2 $\,\,\,W(\eta) \cdot \varphi =0\,$ for each 2-form $\eta$. 
Using (\ref{5.2}) we obtain
\[
\left( \lambda_j + \textstyle{\frac{1}{n-2}}\, \Lambda_0 \right) \fa_j \cdot 
V \cdot \varphi =0 \, \,\,\quad \quad j=1, \ldots, n-2 .
 \]
Since $M_{\lu}$ is not conformally flat and $\fa_j$ is space-like, it
follows
\bea \label{5.9}
0= V \cdot \varphi
\eea
Furthermore, we have
\[
\nabla^S_X D \varphi \quad \stackrel{(\ref{4.4})}{=} \,
\textstyle{\frac{n}{2}}\, K(X) \cdot \varphi .
 \]
Using (\ref{5.1}) and (\ref{5.9}) we obtain
\[
\nabla^S_X D \varphi = - \textstyle{\frac{n}{2(n-2)}\, \Lambda_0 
\,g_{\lu}(X,V)V 
\cdot \varphi =0 }.
 \]
Hence, $\,D \varphi\,$ is parallel. From Proposition 5 follows, that $D 
\varphi$ is constant and $\,V \cdot D \varphi =0\,$. Then the twistor
equation and (\ref{5.3}) yield
\[
 \textstyle{0= \nabla^S_V \varphi + \frac{1}{n} V \cdot D \varphi =
\frac{\partial}{\partial t} (\varphi)} .
 \]
Therefore, $\varphi$ does not depend on $t$. Moreover, the twistor equation 
implies that
\[
 - \fa_{\overline{0}} \cdot \nabla^S_{\fa_{\overline{0}}}
\varphi = \fa_0 \cdot \nabla^S_{\fa_0} \varphi = \frac{1}{n} D
\varphi .
 \]
Then the formulas (\ref{5.6}), (\ref{5.7}) and (\ref{5.9}) give
\[
D \varphi = - n \textstyle{\fa_{\overline{0}} \cdot \left(
\frac{\partial}{\partial s} - \frac{1}{2} \Big(
\sum\limits^{n-2}_{j=1} \lambda_j x^2_j +1 \Big)
\frac{\partial}{\partial t} \right) (\varphi)}
 = n \textstyle{\fa_{0} \cdot \left( \frac{\partial}{\partial s} -
\frac{1}{2} \Big( \sum\limits^{n-2}_{j=1} \lambda_j x^2_j -1
\Big) \frac{\partial}{\partial t} \right) (\varphi)}.
\]
Since $\varphi$ does not depend on $t$ we obtain
\[
2D \varphi = \textstyle{n (\fa_0 - \fa_{\overline{0}}) \cdot
\frac{\partial}{\partial s} (\varphi) \,
=\,n V \cdot \frac{\partial}{\partial s} (\varphi)\, =\, n
\frac{\partial}{\partial s} (V \cdot \varphi ) \stackrel{(\ref{5.9})}{=} 0}. 
\]
Therefore, $\varphi$ is harmonic and the twistor equation implies
that $\varphi$ is parallel.
\qed
\vspace{0.5cm}\\
Now, let $(M^n,g)$ be a non-conformally flat, non-simply connected,
indecomposable Lorentzian symmetric space of dimension $n \ge 3$.
Then, according to the Theorems 1,2 and 4, $(M^n,g)$ is isometric to
$\,M^n_{\lu}/A\,$, where $A$ is a discrete subgroup of the centralizer
$Z_{\lu} := Z_{I(M_{\lu})} (G(M_{\lu}))$.\\

{\bf 1. case:}
There exist $i \in \{1, \ldots, n-2 \}$ such
that $\lambda_i >0$ or $(i,j)$ such that
$\frac{\lambda_i}{\lambda_j} \not\in {\Bbb Q}^2$. Then $Z_{\lu} \simeq
{\Bbb R} = \left\{ \gamma_{\alpha} \mid \gamma_{\alpha} (s,t,x)=(s,t + \alpha , 
x),
\,\, \alpha \in {\Bbb R} \right\}$. \\
Let $\gamma \in A$. With respect to the global basic
$(\fa_{\overline{0}}, \fa_0, \fa_1, \ldots, \fa_{n-2})$ the
differential $d \gamma_y$ corresponds to the matrix $\Gamma (y)
\equiv E \in SO(n,1)$. Hence, $\tilde{\Gamma}^{\pm} (y)= \pm 1 \in
Spin (n,1)$. Therefore, we have 2 spin structures on
$\,M=M_{\lu}/A\,$ corresponding to the homomorphisms $\mbox{Hom} (A;{\Bbb 
Z}_2)$. If $\,\varepsilon \in \mbox{Hom} (A, {\Bbb Z}_2)\,$ is not trivial, there 
are no non-trivial $\varepsilon$-invariant constant spinor fields on 
$M_{\lu}$. 
From the Propositions 3,5 and 6 it follows that the twistor spinors
on $M=M_{\lu}/A$ are given by
\[
{\cal T} (M_{\lu}/A, \varepsilon)\,= \,\left\{ \begin{array}{ll} \{
\varphi \in C^{\infty} (M, \Delta_V) \mid \, \varphi \,\, \mbox{const} \} 
& \,\varepsilon \, \mbox{trivial}\\[0.1cm] \{0\} & \,\varepsilon \,
\mbox{non-trivial}. \end{array} \right.
 \]
{\bf 2. case:}
Let $\lambda_j = - k_i^2 <0$ and
$\,\frac{k_i}{k_j} \in {\Bbb Q}\,$ for all $i,j= 1, \ldots, n-2.$
Then
\begin{eqnarray*} Z_{\lu}& \simeq & \left\{ \begin{array}{rl}
\gamma_{\underline{m}, \alpha} \Big| & \gamma_{\underline{m},
\alpha} (s,t,x)\,:=\, (s+ \beta, t+ \alpha\,, (-1)^{m_1} x_1, \ldots,
(-1)^{m_{n-2}} x_{n-2}); \\[0.1cm]
\Big| & \,\mbox{where} \, \, \alpha \in {\Bbb R}, \,\underline{m} =(m_1, 
\ldots, m_{n-2}) \in {\Bbb Z}^{n-2},\\
\Big| & \,\beta = \pi \cdot \frac{m_i}{k_i} \, \, \, i=
1, \ldots , n-2 \end{array}\right\}\\
\mbox{}\\
& \simeq & {\Bbb Z} \oplus {\Bbb R}.
\end{eqnarray*}
A discrete subgroup $A_{\underline{m}, \alpha} \subset Z_{\lu}$ is
generated by $\gamma_{\underline{m},0} $ and $\gamma_{0, \alpha}$. Let us 
suppose that $\,\sum\limits^{n-2}_{i=1} m_i\,$ is even since otherwise
$\,M_{\lu}/A_{\underline{m}, \alpha}\,$ is not orientable.
 $\,(d \gamma_{\underline{m},
\alpha})_y$ corresponds to the matrix
\[ \Gamma (y)= \left( \begin{array}{ccccc} 1&&&&0 \\
                                                & 1 &&&\\
                                                && (-1)^{m_1} && \\
                                                &&& \ddots &\\
                                                0 &&&& (-1)^{m_{n-2}}
\end{array} \right) .
 \]

Hence $\tilde{\Gamma}^{\pm} (y)= \pm e_1^{m_1} \cdot \ldots \cdot
e^{m_{n-2}}_{n-2}$. Let $\,m_{i_1}, \ldots , m_{i_s}\,$ be the odd
elements in the tupel $\underline{m}\,\, (s \in 2 {\Bbb Z})$, and let us 
denote by $\,\omega_{\underline{m}} \in Spin(n-2,0) \subset Spin(n,1)\,$ the 
element 
\[
\omega_{\underline{m}} = e_{i_1}  \cdot \ldots \cdot e_{i_s} .
 \]
Then because of $\,\omega^2_{\underline{m}} =(-1)^{\frac{s}{2}}\,$,
$\omega_{\underline{m}}$ is an involution on $\Delta_{n-2,0}$ if $s
\equiv 0(4)$ and an almost complex structure if $s \equiv 2(4)$.
The eigenspaces of $\omega_{\underline{m}}$ to the eigenvalues $\pm 1$
and $\pm i$, respectively, have the same dimension (see formulas
(\ref{phi1}) and (\ref{phi2})).\\
The manifold $\,M_{\lu}/A_{\underline{0}, \alpha} , \,\alpha \not= 0\,$, has 2 
spin structures and the twistor spinors are given as in case 1.
$\,M_{\lu}/A_{\underline{m}, 0} , \,\underline{m} \not=
\underline{0}\,$, has 2 spin structures, described by the
homomorphisms
 $\,\varepsilon_{\pm} \in \mbox{Hom}  \, (A_{\underline{m},0} , Spin
(n-2,0) )\,$ given by
$\,\varepsilon_{\pm} (\gamma_{\underline{m},0} )= \pm
\omega_{\underline{m}}\,$ .
Then, according to the Propositions 3,5 and 6 the twistor spinors
on $\,M=M_{\lu}/A_{\underline{m},0}\,$ are
\[
 {\cal T} (M_{\lu}/A_{\underline{m},0} , \varepsilon_{\pm} )\,=\,
\left\{ \begin{array}{cl} \varphi \in C^{\infty} (M_{\lu}, \Delta_{n,1}) 
\, \,  \mid & \,\varphi (y) \equiv v \otimes u(-1);\\[0.1cm]
&  \mbox{where} \, v \in \Delta_{n-2,0} \,\,\mbox{and}\,\,
\omega_{\underline{m}} \cdot v = \pm v \end{array}
\right\} .
 \]
Hence,
\[
\dim {\cal T}  (M_{\lu}/A_{\underline{m},0} \, ,
\varepsilon_{\pm})\,= \,\left\{ \begin{array}{ll} \,0 & \, \mbox{if}\,\, s 
\equiv 2(4)\\[0.1cm]
\, \textstyle{\frac{1}{4}} \cdot 2^{\left[\frac{n}{2}\right]} &
\,\mbox{if}\, \,s \equiv 0(4). \end{array} \right.
\]
The manifold $\,M_{\lu}/A_{\underline{m}, \alpha} \, , \,\underline{m} 
\not=0, 
\alpha \not=0\,$, has 4 spin structures corresponding to the homomorphisms 
$\,\var \in \,\mbox{Hom}\,(A_{\underline{m},\alpha};Spin(n,1))\,$ given by 
$\,\var(\gamma_{\underline{m},0}) = \pm \omega_{\underline{m}},\, 
\var(\gamma_{0,\alpha})= \pm 1\,$. Hence, 
\[
{\cal T} (M_{\lu}/A_{\underline{m}, \alpha} \,, \varepsilon )\,= \,
\left\{ \begin{array}{ll} \, \{0\} & \mbox{} \, \, \varepsilon(\gamma_{0, 
\alpha} )= - 1\\
\mbox{}\\
{\cal T} (M_{\lu}/A_{\underline{m},0} \,, \varepsilon_{\pm}) &
\begin{array}{l} \varepsilon (\gamma_{0, \alpha})=1\\
\varepsilon(\gamma_{\underline{m},0}) =\pm \omega_{\underline{m}}.
\end{array} \end{array} \right.
\]
Summing up, we have


\begin{pr}
Let $(M^n,g)$ be an indecomposable,
non-conformally flat Lorentzian symmetric spin manifold of dimension
$n \ge 3$. Then each twistor spinor is parallel and the dimension of
the space of twistor spinors is

\[ \dim {\cal T} (M^n,g)= q \cdot 2^{\left[ \frac{n}{2} \right]} , \]

where $q = 0, \frac{1}{4}$ or $\frac{1}{2}$, depending on the
fundamental group $\pi_1 (M)$ and on the spin structure.
\end{pr}
\mk

\section{Twistor spinors on indecomposable conformally flat
Lorentzian symmetric spaces of non-constant sectional curvature}

According to the Theorems 1, 2 and 4 there are two isometry classes
of indecomposable, conformally flat, simply connected Lorentzian
symmetric spaces of dimension $n \ge 3$ and non-constant sectional
curvature, namly
\[ M^n_{\pm} := ({\Bbb R}^n, g_{\pm}), \]
where $\hspace{3cm} (g_{\pm})_{(s,t,x)} = 2 ds \,dt \pm || x ||^2 ds^2 + 
\sum\limits^{n-2}_{j=1} dx_j^2 $.\\[0.2cm]
Knowing that $\,\dim {\cal T}(M^n_{\pm})= 2 \cdot 2^{\left[\frac{n}{2}
\right]}\,$ we want to describe the twistor spinors explicitly. Let
$\,w_1, w_2, w_3, w_4 \in \Delta_{n-2,0}\,$ and let us denote by
$\,\varphi_{w_1, w_2, w_3, w_4} \in C^{\infty}(M^n_{\pm},
\Delta_{n,1})\,$ the following smooth functions
\begin{eqnarray*}
\varphi_{w_1, w_2, w_3, w_4} (s,t,x)&:=& (\mp f'(s)w_3 - f(s)w_4 +
x \cdot w_1) \otimes u(1)\\
&& + [-2 w_1 t+ w_2+ x \cdot (f(s)w_3 + f'(s) w_4)] \otimes u(-1),
\end{eqnarray*}
where
\[ f(s)= \left\{ \begin{array}{ll} \mbox{sinh} (s) & \,\mbox{for}\, \,
M^n_{+}\\
\mbox{sin} (s) & \,\mbox{for} \,\, M^n_{-}.
\end{array}  \right. \]

\mk
\begin{pr} The twistor spinors on $\,M^n_{\pm}\,$ are
\[
 {\cal T} (M^n_{\pm}) = \{ \varphi_{w_1, w_2, w_3, w_4} \,\mid \,
w_1, w_2, w_3, w_4 \in \Delta_{n-2,0} \} .
 \]
\end{pr}
\underline{Proof:}  We use the identification
\[
\begin{array}{ccc}
\Delta_{n,1} \simeq \Delta_{n-2,0} \otimes \Delta_{2,1} &
\stackrel{\sim}{\longrightarrow} &
\Delta_{n-2,0} \oplus \Delta_{n-2,0} \\[0.1cm]
 \varphi = \varphi_1 \otimes u(1) + \varphi_2 \otimes u(-1) &
\longmapsto & (\varphi_1, \varphi_2)  .
\end{array}
\]
Then according to (\ref{phi1}) and (\ref{phi2}) the Clifford multiplication 
corresponds to\\

\hspace*{1cm} $ X \cdot \varphi =(-X \cdot \varphi_1, X \cdot \varphi_2) \quad 
\quad \mbox{if $X \in \,\,\mbox{span}\, (\fa_1, \ldots,
\fa_{n-2})$}$,\\[0.2cm]
\hspace*{1cm} $\fa_{\overline{0}} \cdot \varphi =( - \varphi_2, - \varphi_1) 
\quad , \quad \,
\fa_0 \cdot \varphi =( - \varphi_2, \varphi_1) \quad , \quad
V \cdot \varphi =(0, 2 \varphi_1)$.\\

For the spinor derivative we obtain\\

\hspace*{1cm} $\nabla_{\fa_{\overline{0}}}^S \varphi =( \fa_{\overline{0}} 
(\varphi_1), \fa_{\overline{0}} (\varphi_2 ) \pm x \cdot
\varphi_1) \quad , \quad
\nabla^S_{\fa_0} \varphi =(\fa_0 (\varphi_1),
\fa_0(\varphi_2) \pm x \cdot \varphi_1)$,\\[0.2cm]
\hspace*{1cm} $\nabla^S_{\fa_k} \varphi =(\fa_k (\varphi_1), \fa_k
(\varphi_2)) \quad ,\quad k=1, \ldots, n-2$.\\

Let $\varphi=(\varphi_1, \varphi_2)$ be a twistor spinor on
$M^n_{\pm}$. Then, according to Proposition 1, there exists a
spinor field $\psi=(\psi_1, \psi_2)$ on $M^n_{\pm}$ such that
\bea \label{5.10}
(\psi_1, \psi_2)= \fa_k \cdot \nabla^S_{\fa_k}
\varphi =( - \fa_k \cdot \fa_k (\varphi_1), \fa_k \cdot
\fa_k(\varphi_2))
\eea
for each $k=1,2, \ldots, n-2\,$.
Therefore, $\,\varphi_1 (s,t, \cdot)\,$ and $\,\varphi_2 (s,t, \cdot)\,$ are 
twistor spinors on the Euclidean space $\R^{n-2}$. Hence,
\[
 \varphi_i (s,t,x) = u_i (s,t) + x \cdot v_i(s,t) \quad , \quad
 i=1,2 \,,
\]
where $u_i, v_i : {\Bbb R}^2 \to \Delta_{n-2,0}\,$.
From (\ref{5.10}) follows
\[
\psi_1 (s,t,x) = v_1 (s,t) \quad \mbox{and}\quad \psi_2 (s,t,x) = -
v_2 (s,t) .
\]
Furthermore, $\,\psi =( \psi_1, \psi_2)\,$ satisfies
\[
(\psi_1, \psi_2)= - \fa_{\overline{0}} \cdot
\nabla^S_{\fa_{\overline{0}}}  \varphi = \fa_0 \cdot
\nabla^S_{\fa_0} \varphi .
 \]
Therefore,
\bea \label{5.11}
 v_1 & =& \textstyle{\left(
\frac{\partial}{\partial s} - \frac{1}{2} \left( \pm || x||^2 +1
\right) \frac{\partial}{\partial t} \right) (u_2 + x \cdot v_2)\, \pm \, 
x \cdot (u_1 + x \cdot v_1)} \\
\label{5.12}
 v_1 & = & \textstyle{- \left(
\frac{\partial}{\partial s} - \frac{1}{2} \left( \pm ||x||^2 - 1
\right) \frac{\partial}{\partial t} \right) (u_2 + x \cdot v_2) \, \mp \, 
x \cdot (u_1 + x \cdot v_1)}\\
\label{5.13}
- v_2 & = & \textstyle{\left(
\frac{\partial}{\partial s} - \frac{1}{2} \left( \pm ||x||^2 +1
\right) \frac{\partial}{\partial t} \right) (u_1 + x \cdot v_1) }
\\
\label{5.14}
 - v_2 &= & \textstyle{\left(
\frac{\partial}{\partial s} - \frac{1}{2} \left( \pm ||x||^2 - 1
\right) \frac{\partial}{\partial t} \right) (u_1 + x \cdot v_1) }.
\eea
(\ref{5.11}) + (\ref{5.12}) gives
$\,\,\, 2 v_1 = - \frac{\partial}{\partial t} u_2 - x \cdot
\frac{\partial}{\partial t} v_2 \,\,$
and after differentation with respect to $x_k$
$\,\,\, 0= - \fa_k \cdot \frac{\partial}{\partial t} v_2 \,\,$.
Hence,
\bea \label{5.15}
\textstyle{\frac{\partial}{\partial t} v_2 = 0
\quad \mbox{and} \quad v_1 = - \frac{1}{2} \frac{\partial}{\partial
t} u_2}.
\eea
Using this, we obtain from (\ref{5.12}) - (\ref{5.11})
\bea
\label{5.16}
\textstyle{ 0= \mp x \cdot u_1 -
\frac{\partial}{\partial s} u_2 - x \cdot \frac{\partial}{\partial
s} v_2}.
\eea
Differentiation and (\ref{5.15}) show that
\bea
\label{5.17}
\textstyle{ u_1 = \mp \frac{\partial}{\partial s} v_2 \quad
\mbox{and} \quad \frac{\partial}{\partial t}  u_1 =0 }.
\eea
Inserting this in (\ref{5.16}) and using (\ref{5.15}) we obtain
\bea
\label{5.18}
\textstyle{ \frac{\partial}{\partial s} u_2 =0 \quad \mbox{and}
\quad \frac{\partial}{\partial s} v_1 = 0}.
\eea
Hence, $\,u_2 = u_2 (t), \,v_1 = v_1 (t), \,v_2 = v_2 (s)\,$ and $
\,u_1 = u_1(s)\,$. (\ref{5.14}) - (\ref{5.13}) shows
\[
\textstyle{ 0= x \cdot \frac{\partial}{\partial t} v_1 +
\frac{\partial}{\partial t} u_1 = x \cdot v_1' (t)} .
 \]
Therefore, we have $\,v_1 (t) \equiv w_1 \in \Delta_{n-2}\, $ and, because of 
(\ref{5.15}), $\,u_2 (t) = - 2tw_1 + w_2\,$. \\
(\ref{5.14}) + (\ref{5.13}) yields
\[
 - 2 v_2 = 2 u'_1 (s) \mp ||x||^2 \,x \cdot v'_1 (t) = 2 u'_1 (s) ,
 \]
so that, regarding (\ref{5.17}), $\,v_2 (s)= \pm v''_2 (s)\,$. Therefore, 
\[
 v_2 (s)= f(s) w_3 + f'(s) w_4 \qquad \mbox{and} \qquad u_1 (s) =
\mp f'(s) w_3 - f(s) w_4 ,
\]
where
\[
f(s) = \left\{ \begin{array}{ll} \mbox{sinh} (s) & \,\mbox{for} \,\,
\, M^n_{+}\\ \sin (s) & \,\mbox{for} \,\, \, M^n_{-}. \end{array} \right. 
\]
Consequently, the twistor spinor $\varphi$ is of the form $\,\varphi
= \varphi_{w_1, w_2, w_3, w_4}$ . \qed
\vspace{0.7cm}\\
Now, let $(M^n,g)$ be an indecomposable, conformally flat
non-simply connected Lorentzian symmetric space of dimension $n \ge
3$ and non-constant sectional curvature. Then $(M^n,g)$ is
isometric to  $M^n_{+}/A$ or to $M^n_{-}/A$, where $A$ is a discrete
subgroup of
\[
 Z_{+}:=Z_{I(M_{+})} (G (M_{+}))= \{ \varphi_{\alpha} \mid \, \,
\varphi_{\alpha} (s,t,x) = (s,t + \alpha,x)\,; \quad \alpha \in
{\Bbb R}\}
 \]
in the first and of
\[
Z_{-} :=Z_{I(M_{-})} (G(M_{-}))= \left\{ \begin{array}{ll} \varphi_{m,
\alpha} \,\mid & \,\varphi_{m , \alpha} (s,t,x)=(s+ m \pi, t+ \alpha,
(-1)^m x);\\[0.1cm]
&  m \in {\Bbb Z},\,
\alpha \in {\Bbb R}\end{array} \right\}
 \]
in the second case.\\

{\bf 1. case:} $\,M=M^n_{+}/A_{\alpha} \,,\, A_{\alpha} = {\Bbb Z}
\varphi_{\alpha}\,$.
Then there are 2 spin structures corresponding to  $\,\varepsilon
\in \mbox{Hom} (A_{\alpha} , {\Bbb Z}_2)\,$.  The Propositions 3
and 8 show
\[
 {\cal T} (M, \varepsilon)= \left\{ \begin{array}{ll} \{
\varphi_{0, w_2, w_3, w_4} \,\mid
\,w_2, w_3, w_4 \in \Delta_{n-2,0} \} & \,\,
\varepsilon =1\\[0.1cm]
\{0\} & \,\,\varepsilon \not= 1. \end{array} \right.
 \]
{\bf 2. case:} $\,M=M^n_{-}/A_{m, \alpha} , \,A_{m, \alpha} =
\langle \varphi_{m,0} , \varphi_{0, \alpha} \rangle\, $.
If $m$ is even and $\alpha \not= 0$, we have the same result as in
case 1, since $f(s)= \sin (s)$ is $2 \pi {\Bbb Z}$-invariant. If $m$
is odd, $M$ is orientable only if $n$ is even. Then $M^{2k}$ has 2
spin structures if $\alpha =0$ and 4 spin structures if $\alpha
\not=0$. The Propositions 3 and 8 show
\[
 {\cal T} (M^{2k}_{-}/A_{m,0} , \varepsilon)= \left\{
\begin{array}{ll} \{0\} & \,n=2k \equiv 0(4) \\[0.3cm]
\left\{ \begin{array}{l} \varphi_{w_1, w_2, w_3, w_4} \,\mid \\
\qquad w_1, w_2 \in \Delta^{\pm}_{n-2,0} \\
 \qquad w_3, w_4 \in \Delta^{\mp}_{n-2,0} \end{array} \right\} &
\begin{array}{l} n=2k \equiv 2(4)\\ \varepsilon (\varphi_{m,0})=
\pm e_1 \cdot \ldots \cdot e_{n-2} \end{array}  \end{array} \right.
\]
\mk
\[
 {\cal T} (M^{2k}_{-}/A_{m,\alpha} , \varepsilon)= \left\{
\begin{array}{ll} \{0\} & \begin{array}{l} \varepsilon
(\varphi_{0,\alpha})=-1 \\
\mbox{or} \, \, n=2k \equiv 0(4) \end{array} \\[0.4cm]
 \left\{\begin{array}{l} \varphi_{0, w_2, w_3,  w_4} \mid  \\
\qquad w_2 \in \Delta^{\pm}_{n-2,0} \\
\qquad w_3, w_4 \in \Delta^{\mp}_{n-2,0} \end{array} \right\} &
\begin{array}{l} n=2k \equiv 2(4)\\ \varepsilon (\varphi_{0,
\alpha}) = 1\,\,\mbox{and}\\
\varepsilon (\varphi_{m,0}) = \pm e_1 \cdot \ldots \cdot e_{n-2}.
\end{array}  \end{array} \right.
 \]
\vspace{0.5cm}\\
Summing up, we have in particular:

\begin{pr}  Let $(M^n,g)$ be an indecomposable,
conformally flat Lorentzian  spin manifold $(M^n,g)$ of non-constant
sectional curvature and dimension $n \ge 3$. Then the dimension of
the space of twistor spinors is
\[
\dim \cT(M^n,g)= q \cdot 2^{\left[ \frac{n}{2} \right] } ,
 \]
where $q =0, \frac{3}{4}, 1, \frac{3}{2}$ or $2$, depending on the
fundamental group $\pi_1 (M)$ and on the spin structure.
\end{pr}
\mk


\section{Twistor spinors on Lorentzian symmetric spaces of constant
sectional curvature}

Let $\psi_{u,v} \in C^{\infty} ({\Bbb R}^{n+1,k}, \Delta_{n+1,k})$
denote the twistor spinors on the pseudo-Euclidean space ${\Bbb
R}^{n+1,k}$
\[
\psi_{u,v} (x):= u+x \cdot v \quad , \quad u,v \in \Delta_{n+1,k} .
\]
The pseudo-sphere $S^n_1 (r) \subset {\Bbb R}^{n +1,1}$ and the
pseudo-hyperbolic space $H^n_1 (r)$ are umbilic hypersurfaces. Using
the identification of the spinor bundle of the hypersurface with
that of the external space (see section 4) we obtain from
Proposition 4:

\begin{pr}
The twistor spinors on $S^n_1 (r)$ and $H^n_1 (r)$ with the induced
spin structure are
\bean
{\cal T}(S^n_1 (r)) & = & \left\{ \psi_{u,v}|_{S^n_1(r)} \,\, \Big|
 \,\,
\begin{array}{ll}
 u,v \in \Delta_{n+1,1}  & \,
\mbox{if } \, \, n \equiv 0(2)\\[0.1cm]
 u \in \Delta^+_{n+1,1} ,\, v \in \Delta^-_{n+1,1} & \,
\mbox{if } \, \, n \equiv 1(2) \end{array} \right\} \\[0.3cm]
{\cal T}(H^n_1 (r)) & = & \left\{ \psi_{u,v}|_{H^n_1(r)} \,\, \Big|
 \,\,
\begin{array}{ll} u,v \in \Delta_{n+1,2}
& \, \mbox{if } \,
\, n \equiv 0(2)\\[0.1cm]
  u \in \Delta^+_{n+1,2} ,\, v \in \Delta^-_{n+1,2} & \,
\mbox{if } \, \, n \equiv 1(2) \end{array} \right\} .
\eean
\end{pr}

The Lorentzian manifold $S^n_1 (r)/_{\{\pm I\}}$ is orientable if and
only if $n$ is odd, hence let $n$ be odd. The volume form $\,\omega_{n+1,1} = 
e_1 \cdot \ldots \cdot e_{n+1} \in Spin (n+1,1)\,$ satisfies $
\,\omega^2_{n+1,1}= (-1)^{\frac{n+1}{2} +1}\,$. Therefore,
$\,S^n_1 (r)/_{\{\pm I\}}\,$ has no spin structure, if $n \equiv 3(4)$
and 2 spin structures, if $n \equiv 1(4)$. The spinor fields to these
different spin structures can be identified with the invariant
functions $C^{\infty} (S^n_1 (r),
\Delta^+_{n+1,1})^{\varepsilon_{\pm}}$, where $\varepsilon_{\pm}$ is
the ${\Bbb Z}_2$-action given by
\[
(\varepsilon_{\pm} (-1) \varphi)(x)= \pm \omega_{n+1,1} \cdot \varphi
(-x)= \pm \varphi (-x) .
\]
From the Propositions 3 and 10 follows
\[
{\cal T} (S^{4k+1}_1 (r) /_{\{\pm I\}} , \varepsilon )\,=\, \left\{
\begin{array}{cc} \big\{ \psi_{u^+,0}|_{S^{4k+1}_1}\,\, \,\big|\,\, \, u^+ 
\in \Delta^+_{4k+2,1} \big\} & 
\,\,\mbox{if  } \varepsilon = \varepsilon_+ \\[0.2cm]
\big\{ \psi_{0,v^-}|_{S^{4k+1}_1}\, \,\,\big|\,\, \, v^- \in
\Delta^-_{4k+2,1} \big\} & \,\,\mbox{if } \varepsilon = \varepsilon_-.
\end{array}
\right.
\]
Now, let us consider a Lorentzian symmetric space $M^n$ of constant
negative sectional curvature. Then $M^n$ is isometric to a
Lorentzian covering of $H^n_1 (r) /_{\{\pm I\}}$. Let
\[
\begin{array}{lccc}
\tilde{\pi} :& \widetilde{H^n_1 (r)} = {\Bbb R} \times {\Bbb
R}^{n-1} & \longrightarrow  & H^n_1 (r) \subset {\Bbb R}^{2,2} \times
{\Bbb R}^{n-1} \\[0.2cm]
& (t,x) & \longmapsto & \textstyle{(\sqrt{r^2+||x||^2} \cos{t} ,
\sqrt{r^2+||x||^2} \sin{t} , x)}
\end{array}
\]
be the universal Lorentzian covering of $H^n_1 (r)$.
Let $\hat{Q}$ denote the reduction of the trivial spin structure
$Q$ of ${\Bbb R}^{n+1,2}$ to the subgroup $Spin (n,1)$ given by the
Gaus map. Then $\tilde{Q} := \tilde{\pi}^* \hat{Q}$ is the uniquely
determined spin structure of $\widetilde{H^n_1} (r)$. The spinor
fields of $\widetilde{H^n_1} (r)$ can be identified with the smooth
functions $\,C^{\infty} ( \widetilde{H^n_1} (r),
\Delta^{(\pm)}_{n+1,2})\,$, the twistor spinors are given by
\[
{\cal T} (\widetilde{H^n_1} (r)) = \{ \tilde{\psi}_{u,v} :=
\psi_{u,v}|_{H^n_1(r)} \circ \tilde{\pi} \,\, \Big| \,\,\, 
\psi_{u,v}|_{H^n_1(r)} \in \cT(H^n_1 (r)) 
\} .
\]
Let
\[
\begin{array}{lccc}
\pi_m :& N^n_m & \longrightarrow & H^n_1 (r) \\[0.2cm]
& \textstyle{( \sqrt{. }\, \cos{t} , \sqrt{. }\, \sin{t} , x)} & \longmapsto & 
 \textstyle{(\sqrt{. }\, \cos{(mt)}, \sqrt{. }\, \sin{(mt)}, x)},
\end{array}
\]
($\sqrt{. }=\sqrt{r^2+||x||^2})\,\,$, be the Lorentzian covering of
$H^n_1 (r)$ with respect to $m {\Bbb Z} \subset \pi_1 (H^n_1 (r)) =
{\Bbb Z},\, m=1,2,3, \ldots.\, $.
The manifold $N^n_m$ has 2 spin structures. The corresponding spinor
fields are
given by the $\varepsilon^{\pm}_m$-invariant  functions
$\,C^{\infty} (\widetilde{H^n_1 (r)} , \Delta^{(\pm)}_{n+1,2}
)^{\varepsilon^{\pm}_m}\, $, where $\varepsilon^{\pm}_m$ is the $m
{\Bbb Z} $-action
\[
(\varepsilon^{\pm}_m (mz) \varphi)(t,x)=(\pm 1)^z \varphi (t+ 2
\pi mz, x) .
\]
Therefore, the twistor spinors on $N^n_m$ are
\[
{\cal T} (N^n_m, \varepsilon)\,=\, \left\{ \begin{array}{ll} \{
\psi_{u,v}|_{H^n_1(r)} \circ \pi_m \,\,\, \Big| \,\, \psi_{u,v}|_{H^n_1(r)} 
\in \cT(H^n_1 (r)) \} & 
\,\,\mbox{if }\,
\varepsilon = \varepsilon^+_m \\[0.1cm]
\{0\} &\,\,\mbox{if }\, \varepsilon = \varepsilon^-_m.
\end{array} \right.
\]
\mk
Finally, let us consider the manifolds $\, N^n_m /_{\{\pm I\}}\,$.
Since $\,N^n_m /_{\{\pm I\}}\,$ is orientable if and only if $n$ is
odd, let $n$ be odd. For the volume form $\,\omega_{n+1,2} = e_1 \cdot
\ldots \cdot e_{n+1} \in Spin (n+1,2)\,$ we have  $\,\omega^2_{n+1,2}
=(-1)^{\frac{n+1}{2} +2}\,$. Therefore, there is no spin
structure on $N^n_m /_{\{\pm I\}}$ if  $n \equiv 1(4)$ and there are
4 spin structures in case $n \equiv 3(4)$. The spinor fields are
given by the functions $\,C^{\infty} (\widetilde{H^n_1 (r)},
\Delta^+_{n+1,2})^{(\varepsilon^{\pm}_m, \delta^{\pm})}\,$, invariant
under the $m {\Bbb Z}$-action $\varepsilon^{\pm}_m$ and the ${\Bbb
Z}_2$-action $\delta^{\pm}$ given by
\begin{eqnarray*}
(\delta^{\pm} (-1) \varphi)(t,x)&=& \pm \omega_{n+1,2} \cdot \varphi
(t+m \pi , -x) \\
&=& \pm \varphi (t + m \pi , -x ).
\end{eqnarray*}
Then, the twistor spinors are
\[
{\cal T}  \left(N^{4k+3}_m /_{\{ \pm I\}} , \varepsilon \right)
\,=\, \left\{ \begin{array}{ll} \{0\} & \begin{array}{l} \varepsilon =( 
\varepsilon^-_m , \delta^{\pm}) \,\,\mbox{or}\\ \varepsilon = ( 
\varepsilon^+_m, 
\delta^-),\, m \equiv 0(2) \end{array} \\[0.4cm]
\{ \tilde{\psi}_{u^+,0} \,\mid\, u^+ \in \Delta^+_{4k+4,2} \} & \mbox{} \, 
\, \varepsilon = ( \varepsilon^+_m , \delta^+ )\\[0.3cm]
\{ \tilde{\psi}_{0,v^-} \,\mid \, v^- \in \Delta^+_{4k+4,2} \} &
\,\varepsilon =( \varepsilon^+_m , \delta^-),\, m \equiv 1(2).
 \end{array} \right.
 \]
Summing up, we have in particular

\begin{pr} Let $(M^n,g)$ be a Lorentzian symmetric
spin manifold of constant sectional curvature $k \not=0$ and dimension
$n \ge 3$, then the dimension of the space of twistor spinors is
\[
\dim \cT (M^n,g)= q \cdot 2^{\left[ \frac{n}{2} \right]}  ,
\]
where $q =0,1,$ or 2 depending on $\pi_1 (M)$ and on the spin
structure.
\end{pr}
\mk

\bibliographystyle{alpha}
\bibliography{lit}
\vspace{1cm}
{\footnotesize Helga Baum\\
Institut f\"ur Reine Mathematik\\
Humboldt-Universit\"at zu Berlin\\
Sitz: Ziegelstr. 13a\\
10099 Berlin , Germany\\
email: baum@mathematik.hu-berlin.de\\}
\end{document}